\documentclass{article}
\usepackage[utf8]{inputenc}
\usepackage{mathtools} 
\usepackage{amsfonts} 
\usepackage{amsmath}
\usepackage{float}
\usepackage{nccmath}
\title{On the 2-adic valuation of ($a^b-c^d$)}
\author{Luca Onnis}
\date{August 2021}

\begin{document}
\maketitle
\begin{abstract}
    In this paper will be proved an inequality regarding $v_2(a^{b}-c^{d})$. Using this formula it will be possible to have informations about the divisibility of 2 of this function without computing it. Then, will be studied the behavior of this function when it approaches infinity and there will be a lot of analogies with 2-adic integers.
\end{abstract}
\section{Introduction}
There are some fundamental properties about p-adic valuation function reported below. \cite{tdn}
\subsection{Definition 1}
Let $p$ be  a prime number. If $n$ is an integer not equal to 0, $p$-adic valuation of $n$ is equal to:
\[
v_p(n):= \max\Bigl[k\in\mathbb{N}:p^{k}\mid n\Bigr]
\]
\subsection{Lemma 1}
Let $p$ be a prime number and $v_p$ the $p$-adic valuation function. Hence, we'll have that:
\[
v_p(ab)=v_p(a)+v_p(b)
\]
for $a,b\in\mathbb{Z}$ with $a,b$ $\neq 0$; \\
\subsection{Lemma 2} 
Let $p$ be a prime number and $v_p$ the $p$-adic valuation function. Hence it satisfies the following inequality:
\[
v_p(a+b) \geq\min\Bigl[v_p(a);v_p(b)\Bigr]
\]
for $a,b\in\mathbb{Z}$ with $a,b \neq 0$, and $v_p(a)\neq v_p(b)$
\subsection{Lemma 3}
Let $x,y$ be odd integers and $n$ and even integer. Hence:
\[
v_2(x^{n}-y^{n})=v_2(x-y)+v_2(x+y)+v_2(n)-1
\]
\subsection{Lemma 4}
Let $x,y,n$ odd integers. Hence:
\[
v_2(x^{n}-y^{n})=v_2(x-y)
\]
\section{Theorem 1}
For every $a,b,c,d\in\mathbb{N}$ with $a,c$ odd integers, $d$ even and for every integer $b$ we'll have that:
\[
v_2(a^{b}-c^{d}):\begin{cases} \geq\max\Bigl[v_2(a+1);v_2(a-1)\Bigr]+1, & \mbox{if } b\mbox{ is even} \\ = \min\Bigl[v_2(a-1);\max[v_2(c+1);v_2(c-1)]+v_2(d)\Bigr], & \mbox{if } b\mbox{ is odd} \end{cases}
\]
if and only if $\max\Bigl[v_2(c+1);v_2(c-1)\Bigr] + v_2(d) \geq\max\Bigl[v_2(a+1);v_2(a-1)\Bigr]+1$
Consider the following case:
\[
\max\Bigl[v_2(c+1);v_2(c-1)\Bigr]=v_2(c+1)
\]
\subsection{Lemma 5}
If $v_2(c+1)$ is the maximum value between those two, it means that $c+1$ is surely divisible by 4, so $v_2(c+1)\geq 2$ , while $v_2(c-1)=1$. Consider the integers $c-1$, $c$ e $c+1$. For hypothesis c is odd so we can write those integers as $2c-2$, $2c-1$ and $2c$. But $4\mid 2c$ so this implies that $c=2h$ for some $h\in\mathbb{N}$ and so $2c-2=4h-2$ which is not divisible by 4. In fact:
\[
4h-2\equiv 2 \mod(4)
\]
But obviously
\[
4h-2\equiv 0 \mod(2)
\]
Note that if we consider $v_2(c-1)$ as the maximum between those two, $c-1$ will be surely divisible by 4, while $c+1$ will be an even number not divisible by 4.  \\
Furthermore consider:
\[
\max\Bigl[v_2(a+1);v_2(a-1)\Bigr]=v_2(a+1)
\]
and applying the same strategy we'll have that $v_2(a+1)\geq 2$ , while $v_2(a-1)=1$.
So:
\[
v_2(a^{b}-c^{d})=v_2(a^{b}-1+1-c^{d})
\]
For lemma 2 we'll have that:
\[
v_2(a^{b}-1+1-c^{d})\geq\min\Bigl[v_2(a^{b}-1);v_2(1-c^{d})\Bigr]
\]
but $v_2(1-c^{d})=v_2(c^{d}-1)$ So we'll have that:
\[
v_2(a^{b}-c^{d})\geq\min\Bigl[v_2(a^{b}-1);v_2(c^{d}-1)\Bigr]
\]
Because of lemma 3 we'll have that:
\[
v_2(a^{b}-c^{d})\geq\min\Bigl[v_2(a-1)+v_2(a+1)+v_2(b)-1;v_2(c-1)+v_2(c+1)+v_2(d)-1\Bigr]
\]
But $v_2(a-1)=v_2(c-1)=1$ , so the inequality becames:
\[
v_2(a^{b}-c^{d})\geq\min\Bigl[v_2(a+1)+v_2(b);v_2(c+1)+v_2(d)\Bigr]
\]
We can certainly conclude that $v_2(a+1)+v_2(b)\geq v_2(a+1)+1$ because for even $b$ $v_2(b)\geq 1$ \\
Consider the case where $v_2(a+1)+v_2(b)\leq v_2(c+1)+v_2(d)$ ; so:
\[
v_2(a+1)+1\leq v_2(a+1)+v_2(b)\leq v_2(c+1)+v_2(d)
\]
So we'll have that:
\[
v_2(a^{b}-c^{d})\geq v_2(a+1)+1
\]
for every even $b,d$ and every odd $a,c$ such that:
\[
\max\Bigl[v_2(c+1);v_2(c-1)\Bigr] + v_2(d) \geq\max\Bigl[v_2(a+1);v_2(a-1)\Bigr]+1
\]
If $v_2(b)$ is sufficiently large we will have that: 
\[
v_2(a+1)+v_2(b)\geq v_2(c+1)+v_2(d)
\]
hence:
\[
v_2(a+1)+1\leq v_2(c+1)+v_2(d)\leq v_2(a+1)+v_2(b)
\]
so:
\[
v_2(a^{b}-c^{d})\geq v_2(c+1)+v_2(d) 
\]
but $v_2(a+1)+1\leq v_2(c+1)+v_2(d)$ , so surely:
\[
v_2(a^{b}-c^{d})\geq v_2(a+1)+1
\]
Is we consider $v_2(a-1)$ or $v_2(c-1)$ as the maximum the final result would have been:
\[
v_2(a^{b}-c^{d})\geq v_2(a-1)+1
\]
In order to generalize all four combinations below: \\
1) $v_2(a+1)$ max and $v_2(c+1)$ max \\
2) $v_2(a-1)$ max and $v_2(c-1)$ max \\
3) $v_2(a-1)$ max and $v_2(c+1)$ max \\
4) $v_2(a+1)$ max and $v_2(c-1)$ max \\
we'll using the following formula: \\
\[
v_2(a^{b}-c^{d})\geq\max\Bigl[v_2(a+1);v_2(a-1)\Bigr]+1
\]
For every odd $a,c\in\mathbb{N}$ and for every even $b,d\in\mathbb{N}$ if and only if $\max\Bigl[v_2(c+1);v_2(c-1)\Bigr] + v_2(d) \geq\max\Bigl[v_2(a+1);v_2(a-1)\Bigr]+1$ \\
Consider odd b, we'll have that:
\[
v_2(a^{b}-c^{d})\geq\min\Bigl[v_2(a^{b}-1);v_2(c^{d}-1)\Bigr]
\]
Because of lemma 4 we'll obtain that:
\[
v_2(a^{b}-c^{d})\geq\min\Bigl[v_2(a-1);v_2(c+1)+v_2(d)\Bigr]
\]
And generalizing it we'll have that:
\[
v_2(a^{b}-c^{d})\geq\min\Bigl[v_2(a-1);\max[v_2(c+1);v_2(c-1)]+v_2(d)\Bigr]
\]
for every odd integers b. We want to prove that:
\[
v_2(a^{b}-c^{d})=\min\Bigl[v_2(a-1);\max\Bigl[v_2(c+1);v_2(c-1)\Bigr]+v_2(d)\Bigr]
\]
I define:
\[
m=\max\Bigl[v_2(c+1);v_2(c-1)\Bigr]+v_2(d)
\]
So we want to prove that:
\[
v_2(a^{b}-c^{d})=\min\Bigl[v_2(a-1);m\Bigr]
\]
For hypothesis we'll have that:
\[
m\geq\max\Bigl[v_2(a+1);v_2(a-1)\Bigr]+1
\]
So $\min\Bigl[v_2(a-1);m\Bigr]=v_2(a-1)$. In fact if for absurdity $m\leq v_2(a-1)$ we will have that:
\[
\begin{cases}
  m\geq\max\Bigl[v_2(a+1);v_2(a-1)\Bigr]+1 \\
       m\leq v_2(a-1)
\end{cases}
\]
Namely:
\[
\max\Bigl[v_2(a+1);v_2(a-1)\Bigr]+1\leq m \leq v_2(a-1)
\]
which is impossible. So we must have that:
\[
\begin{cases}
  m\geq\max\Bigl[v_2(a+1);v_2(a-1)\Bigr]+1 \\
       m\geq v_2(a-1) \\
       m\geq 2\\
\end{cases}
\]
In fact if c is odd and d is even we'll have that $m\geq 2$. \\
So that $m\geq v_2(a-1)$ is true, $a$ must be in the form:
\[
a=2^{m-j}+1+2^{m-j+1}k
\]
with $k,j\in\mathbb{N}$ and $0\leq j \leq m-2$
We want to prove that for every odd $b$:
\[
v_2(a^{b}-c^{d})=v_2(a-1)
\]
So:
\[
v_2\Biggr(\Bigl[2^{m-j}+1+2^{m-j+1}k\Bigr]^{b}-c^{d}\Biggl)=v_2(a-1)
\]
If this is true:
\[
\begin{cases}
  \Bigl[2^{m-j}+1+2^{m-j+1}k\Bigr]^{b}-c^{d} \equiv 0 \mod(2^{m-j}) \\
       \Bigl[2^{m-j}+1+2^{m-j+1}k\Bigr]^{b}-c^{d} \not\equiv 0 \mod(2^{m-j+1})
\end{cases}
\]
For every odd $b$ we'll have that:
\[
\begin{cases}
  (1-c^{d})\equiv 0 \mod(2^{m-j}) \\
       (2^{m-j}+1-c^{d}) \not\equiv 0 \mod(2^{m-j+1})
\end{cases}
\]
In order to prove that:
\[
 \Bigl[2^{m-j}+1+2^{m-j+1}k\Bigr]^{b}-c^{d} \equiv (2^{m-j}+1-c^{d}) \mod(2^{m-j+1})
\]
for odd $b$, consider Newton's binomial expansion of:
\[
\Bigl[2^{m-j}+1+2^{m-j+1}k\Bigr]^{b} = \Bigl[2^{m-j}(2k+1)+1\Bigr]^{b}
\]
b is odd, so it is possible to write $b$ as $b=2p+1$. Let $m-j=v$ \\
\subsection{Lemma 6}
\[
(a+b)^{n}= \sum_{k=0}^{n}\binom{n}{k}a^{n-k}b^{k}
\]
In our case we want to prove that:
\[
\Bigl[2^{m-j}(2k+1)+1\Bigr]^{b} \equiv \Bigl [2^{v}(2k+1)+1\Bigr]^{2p}\Bigl[2^{v}(2k+1)+1\Bigr] \equiv (2^{v}+1) \mod(2^{v+1})
\]
Which is equal to say that:
\[
\Bigl[2^{v}(2k+1)+1\Bigr]^{2p} \equiv 1 \mod (2^{v+1})
\]
because:
\[
\Bigl[2^{v}(2k+1)+1\Bigr]\equiv 2^{v}+1 \mod (2^{v+1})
\]
Using Newton's formula we'll have that:
\[
\Bigl[2^{v}(2k+1)+1\Bigr]^{2p}= \sum_{i=0}^{2p}\binom{2p}{i}1^{2p-i}\Bigl[2^{v}(2k+1)\Bigr]^{i}=\sum_{i=0}^{2p}\binom{2p}{i}\Bigl[2^{iv}(2k+1)^{i}\Bigr]
\]
So we want to prove that:
\[
\sum_{i=0}^{2p}\binom{2p}{i}\Bigl[2^{iv}(2k+1)^{i}\Bigr] \equiv 1 \mod (2^{v+1})
\]
Note that for $i\geq 2$ , each term is always congruent to 0 mod $2^{v+1}$ because $2^{v+1}|2^{vi}$ and $v+1>vi$ \\
For i=1 the term will be:
\[
\binom{2p}{1}2^{v}(2k+1)=2p*2^{u}(2k+1)=p(2k+1)(2^{v+1}) \equiv 0 \mod (2^{v+1})
\]
For i=0 obiously the term of the sum will be congruent to 1 mod $(2^{v+1})$. So overall we'll have that:
\[
\Bigl [2^{v}(2k+1)+1\Bigr]^{2p} \equiv 1 \mod (2^{v+1})
\]
and so:
\[
\Bigl[2^{m-j}(2k+1)+1\Bigr]^{b} \equiv \Bigl [2^{v}(2k+1)+1\Bigr]^{2p}\Bigl[2^{v}(2k+1)+1\Bigr] \equiv (2^{v}) mod(2^{v+1}) \not\equiv 0 mod (2^{v+1})
\]
for odd b, which is the thesis. \\
We know that $v_2(1-c^{d})= v_2(c^{d}-1)$ and so $(c^{d}-1)  \equiv 0 mod(2^{m-j})$ \\
Because $v_2(c^{d}-1)=m$, so $2^{m}|(c^{d}-1)$ \\ 
And certainly $2^{m-j}|(c^{d}-1)$ \\
But $(c^{d}-1)=2^{m}h$ for some $h\in\mathbb{N}$ and so:
\[
\begin{cases}
(1-c^{d})\equiv 0 mod(2^{m-j}) \\
     (2^{m-j}+2^{m}h) \equiv 2^{m-j} \not\equiv 0 \mod(2^{m-j+1})
\end{cases}
\]
If $j=0$, $m=v_2(a-1)$ and:
\[
v_2(a^{b}-c^{d})=m=v_2(a-1)
\]
For odd $b$.
\section{Theorem 2}
For every $a,b,c,d\in\mathbb{N}$ with odd integers $a,c$, even d and every integer b we'll have that:
\[
v_2(a^{b}-c^{d})=  \begin{cases} \max\Bigl[v_2(c+1);v_2(c-1)\Bigr]+v_2(d), & \mbox{if } b\mbox{ is even} \\ \min\Bigl[v_2(a-1);\max[v_2(c+1);v_2(c-1)]+v_2(d)\Bigr], & \mbox{if } b\mbox{ is odd} \end{cases}
\]
if and only if $\max\Bigl[v_2(c+1);v_2(c-1)\Bigr] + v_2(d) <\max\Bigl[v_2(a+1);v_2(a-1)\Bigr]+1$ \\
As before I define:
\[
\max\Bigl[v_2(c+1);v_2(c-1)\Bigr]=v_2(c+1)
\]
and I consider:
\[
\max\Bigl[v_2(a+1);v_2(a-1)\Bigr]=v_2(a+1)
\]
I define $m=v_2(c+1)+v_2(d)$ \\
For hypothesis $v_2(a+1)+1>m$ , allora:
\[
a=2^{m+h}-1+2^{m+h+1}k
\]
with $h,k\in\mathbb{N}$ \\
We want to prove that for every even integer $b$ we'll have that:
\[
v_2(a^{b}-c^{d})=m
\]
This implies that:
\[
\begin{cases}
  a^{b}-c^{d} \equiv 0 \mod(2^{m}) \\
       a^{b}-c^{d} \not\equiv 0 \mod(2^{m+1})
\end{cases}
\]
But $a=2^{m+h}-1+2^{m+h+1}k$:
\[
\Bigl[2^{m+h}-1+2^{m+h+1}k\Bigr]^{b}-c^{d} \equiv 0 \mod(2^{m})
\]
\[
\Bigl[2^{m+h}(2k+1)-1\Bigr]^{b}-c^{d} \equiv 0 \mod(2^{m})
\]
And for even $b$ we note that, as before:
\[
\Bigl[2^{m+h}(2k+1)-1\Bigr]^{b}-c^{d} \equiv 1-c^{d} \equiv c^{d}-1 \mod(2^{m})
\]
But the following system of congruences must be true:
\[
\begin{cases}
  c^{d}-1 \equiv 0 \mod(2^{m}) \\
     c^{d}-1 \not\equiv 0 \mod(2^{m+1})
\end{cases}
\]
So $v_2(c^{d}-1)=m$ \\
But for lemma 3 we'll have that:
\[
v_2(c-1)+v_2(c+1)+v_2(d)-1=m
\]
But $v_2(c-1)=1$ because of what we proved before. So:
\[
v_2(c+1)+v_2(d)=m
\]
Which is the thesis. \\
So for even integers b we'll have that:
\[
v_2(a^{b}-c^{d})=m=v_2(c+1)+v_2(d)
\]
For odd b we'll again have that:
\[
v_2(a^{b}-c^{d})\geq\min\Bigl[v_2(a^{b}-1);v_2(c^{d}-1)\Bigr]
\]
For lemma 4 we'll have that:
\[
v_2(a^{b}-c^{d})\geq\min\Bigl[v_2(a-1);\max[v_2(c+1);v_2(c-1)]+v_2(d)\Bigr]
\]
for every integer b. But, it will be proved that:
\[
v_2(a^{b}-c^{d})=\min\Bigl[v_2(a-1);\max\Bigl[v_2(c+1);v_2(c-1)\Bigr]+v_2(d)\Bigr]
\]
As the first case consider $\max[v_2(c+1);v_2(c-1)]+v_2(d)$$\leq$ $v_2(a-1)$. For hypothesis we'll have that: \\
$\max[v_2(c+1);v_2(c-1)]+v_2(d)<\max[v_2(a+1);v_2(a-1)]+1$. \\
In this case $v_2(a-1)>v_2(a+1)$. In fact if it was the opposite we would have had that $v_2(a-1)=1$ but because $c$ is odd and $d$ is even $\max[v_2(c+1);v_2(c-1)]+v_2(d)\geq 2$. But this will contradict that $\max[v_2(c+1);v_2(c-1)]+v_2(d)\leq v_2(a-1)$ , we would have had that $2\leq 1$ which is impossible. \\
So $\max[v_2(a+1);v_2(a-1)]=v_2(a-1)$ \\
I define:
\[
m=\max\Bigl[v_2(c+1);v_2(c-1)\Bigr]+v_2(d)
\]
In our case $v_2(a-1) \geq m$ , So $a$ is in the form:
\[
a=2^{m+l}+1+2^{m+l}k
\]
with $l,k\in\mathbb{N}$ \\
As before, if:
\[
v_2(a^{b}-c^{d})=\max\Bigl[v_2(c+1);v_2(c-1)\Bigr]+v_2(d)=m
\]
so:
\[
\begin{cases}
  a^{b}-c^{d} \equiv 0 \mod(2^{m}) \\
       a^{b}-c^{d} \not\equiv 0 \mod(2^{m+1})
\end{cases}
\]
which, as proved before, is true for $a$ in the form $a=2^{m+l}+1+2^{m+l+1}k$. \\
As the second case consider $\max[v_2(c+1);v_2(c-1)]+v_2(d)$$\geq$ $v_2(a-1)$. For hypothesis we'll have that: \\
$\max[v_2(c+1);v_2(c-1)]+v_2(d)<\max[v_2(a+1);v_2(a-1)]+1$. \\
I define:
\[
u=\max[v_2(c+1);v_2(c-1)]+v_2(d)
\]
\[
\begin{cases}
  u\geq v_2(a-1) \\
  u<\max[v_2(a+1);v_2(a-1)]+1
\end{cases}
\]
se $\max[v_2(a+1);v_2(a-1)]=v_2(a+1)$ , allora $v_2(a-1)=1$ e quindi:
\[
\begin{cases}
  u\geq 1 \\
  u< v_2(a+1)+1
\end{cases}
\]
So: 
\[
\min\Bigl[v_2(a-1);u\Bigr]=v_2(a-1)=1
\]
if $\max[v_2(a+1);v_2(a-1)]=v_2(a-1)$ , so $v_2(a+1)=1$ and:\[
\begin{cases}
  u\geq v_2(a-1) \\
  u< v_2(a-1)+1
\end{cases}
\]
Hence:
\[
v_2(a-1)\leq u < v_2(a-1)+1
\]
and so:
\[
u=v_2(a-1)=\max[v_2(c+1);v_2(c-1)]+v_2(d)
\]
and at the end we'll have again that:
\[
\min\Bigl[v_2(a-1);u\Bigr]=v_2(a-1)
\]
So:
\[
v_2(a^{b}-c^{d})=\min\Bigl[v_2(a-1);\max[v_2(c+1);v_2(c-1)]+v_2(d)\Bigr]
\]
For every odd integers $b$. \\
Overall we'll have that: \\
if $\max\Bigl[v_2(c+1);v_2(c-1)\Bigr] + v_2(d) <\max\Bigl[v_2(a+1);v_2(a-1)\Bigr]+1$ \\
\[
v_2(a^{b}-c^{d})=\begin{cases} \max\Bigl[v_2(c+1);v_2(c-1)\Bigr]+v_2(d), & \mbox{if } b\mbox{ is even} \\ \min\Bigl[v_2(a-1);\max[v_2(c+1);v_2(c-1)]+v_2(d)\Bigr], & \mbox{if } b\mbox{ is odd} \end{cases}
\]
and if $\max\Bigl[v_2(c+1);v_2(c-1)\Bigr] + v_2(d) \geq\max\Bigl[v_2(a+1);v_2(a-1)\Bigr]+1$ \\
\[
v_2(a^{b}-c^{d}):\begin{cases} \geq\max\Bigl[v_2(a+1);v_2(a-1)\Bigr]+1, & \mbox{if } b\mbox{ is even} \\ = \min\Bigl[v_2(a-1);\max[v_2(c+1);v_2(c-1)]+v_2(d)\Bigr], & \mbox{if } b\mbox{ is odd} \end{cases}
\]
\section{Examples and considerations}
\subsection{Example 1}
\[
v_2(255^{n}-1023^{2}) \geq 9
\]
for even integers $n$. Because:
\[
\max\Bigl[v_2(1023+1);v_2(1023-1)\Bigr] + v_2(2) \geq\max\Bigl[v_2(255+1);v_2(255-1)\Bigr]+1
\]
\[
10+1\geq 8+1
\]
And for odd integers $n$:
\[
v_2(255^{n}-1023^{2})=1
\]
because
\[
\min\Bigl[v_2(255-1);\max[v_2(1023+1);v_2(1023-1)]+v_2(2)\Bigr]=\min\Bigl[1;11\Bigr]=1
\]
\begin{figure}[h]
    \centering
    \includegraphics[width=8cm]{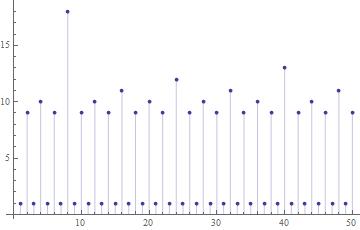}
    \caption{Graph of $v_2(255^{n}-1023^{2})$ from $n$=1 to $n$=50}
    \label{fig:4.1}
\end{figure}
\\
Note that the function is somehow periodic. For example:
\begin{table}[h]
\centering
\begin{tabular}{l | l}
n & $v_2(255^{n}-1023^{2})$ \\
\hline
1+2k & 1 \\
2+4k & 9 \\
4+8k & 10 \\
16k & 11 \\
24+32k & 12 \\
40+64k & 13 
\end{tabular}
\caption{values of $v_2(255^{n}-1023^{2})$ by varying $n$}
\label{tab:4.1}
\end{table} \\
In the fifth section will be explained this behavior, and it will be studied how $v_2(a^{b}-c^{d})$ could diverge to infinity.
\subsection{Example 2}
\[
\begin{cases}
 v_2(6145^{n}-255^{64})\geq 12 & \mbox{for even $n$} \\
       v_2(6145^{n}-255^{64})= 11 & \mbox{for odd $n$}
\end{cases}
\]
In fact we can apply the first theorem, because \[
\max\Bigl[v_2(255+1);v_2(255-1)\Bigr] + v_2(64) \geq\max\Bigl[v_2(6145+1);v_2(6145-1)\Bigr]+1
\]
\[
14 \geq 12
\]
\begin{figure}[h]
    \centering
    \includegraphics[width=8cm]{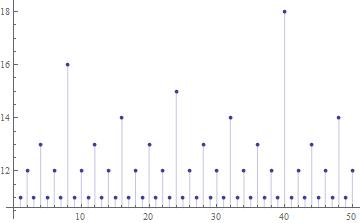}
    \caption{Graph of $v_2(6145^{n}-255^{64})$ from n=1 to n=50}
    \label{fig:4.2}
\end{figure}
\begin{table}[h]
\centering
\begin{tabular}{l | l}
n & $v_2(6145^{n}-255^{64})$ \\
\hline
1+2k & 11 \\
2+4k & 12 \\
4+8k & 13 \\
16k & 14 \\
24+32k & 15 \\
8+64k & 16 
\end{tabular}
\caption{Values of $v_2(6145^{n}-255^{64})$ by varying $n$}
\label{tab:4.1}
\end{table} \\
Note how for different values of $n$ corresponds different values of: 
\[
v_2(6145^{n}-255^{64})
\]
And note how the different values of $n$ are "periodic" in the sense that for different values of $n$ corresponds same values of  the function. This period grows while the function grows, and it is always in the form $T=2^{x}k$
\subsection{Example 3}
\[
\begin{cases}
 v_2(1537^{n}-13^{32}) = 7 & \mbox{for even $n$} \\
       v_2(1537^{n}-13^{32})= 7 & \mbox{for odd $n$}
\end{cases}
\]
In fact we can apply the second theorem, because:
\[
    \max\Bigl[v_2(13+1);v_2(13-1)\Bigr] + v_2(32) <\max\Bigl[v_2(1537+1);v_2(1537-1)\Bigr]+1
\]
\[
7<10
\]
\begin{figure}[h]
    \centering
    \includegraphics[width=8cm]{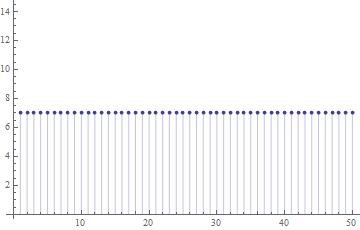}
    \caption{Graph of $v_2(1537^{n}-13^{32})$ from $n$=1 to $n$=50}
    \label{fig:4.3}
\end{figure}
\section{Limits and curiosity}
\subsection{Lemma 7}
It is known that:
\[
v_p(0)=+\infty
\]
From here we can conclude that:
\[
\lim_{v_2(d)\to 0} v_2(a^{b}-c^{d})=v_2(a^{b}-1)
\]
\[
\begin{cases}
 v_2(a^{b}-1) = v_2(a-1)+v_2(a+1)+v_2(b)-1 & \mbox{for even $b$ and odd $a$} \\
       v_2(a^{b}-1) = v_2(a-1) & \mbox{for odd $a,b$}
\end{cases}
\]
but $v_2(d)\to+\infty$ when $d\to 0$ but also when $d$ can be divided infinetely many times by two. \\
Se $d=2^{d_1}k$ , allora $v_2\Bigl(2^{d_1}k\Bigr)\to+\infty$ per $d_1\to+\infty$  \\
In fact I conjecture that if:
\[
v_2(a^{b}-c^{d})=\alpha
\]
\[
v_2(a^{b}-c^{d})=v_2(b)+\max\Bigl[v_2(a+1);v_2(a-1)\Bigr]
\]
for even p and $\alpha\leq v_2(d)+1$ \\
For example, consider the following funtion:
\[
y=v_2(15^{n}-5^{2^{m}})
\]
Because of what I said before $m \to+\infty$ $v_2(15^{n}-5^{2^{m}})\to v_2(n)+4$ \\
\begin{figure}[h]
    \centering
    \includegraphics[width=12cm]{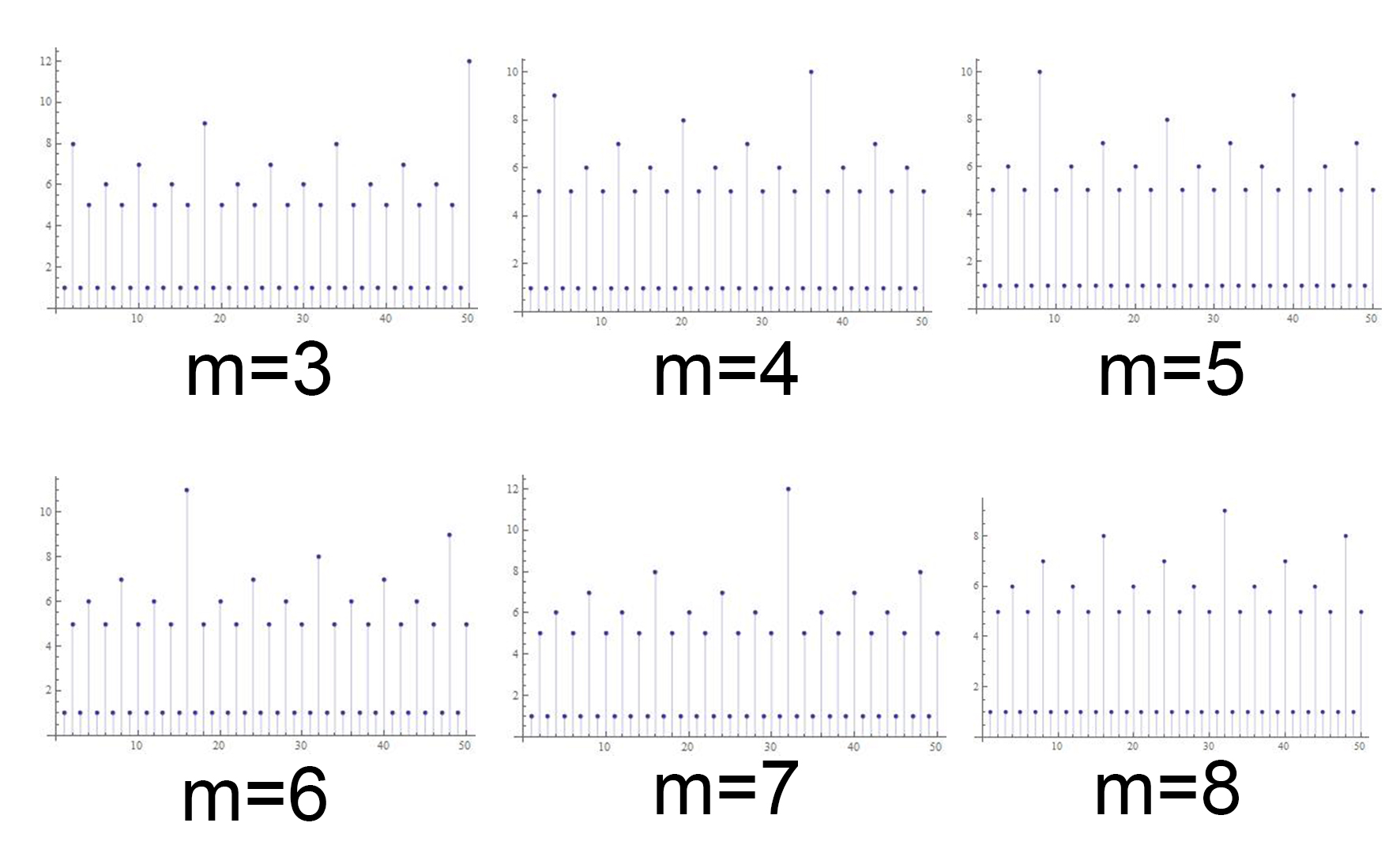}
    \caption{Graph $v_2(15^{n}-5^{2^{m}})$ while $m$ is varying from $m$=3 to $m$=8, from $n$=1 to $n$=50}
    \label{fig:4.3}
\end{figure} \\
\\
\\
\\
\\
\\
\\
\\
\\
\\
\\
\\
\\
\\
Compare the graph of $v_2(15^{n}-5^{2^{8}})$ with the graph of $v_2(n)+4$ reported below: \\
\begin{figure}[H]
    \centering
    \includegraphics[width=7.5cm]{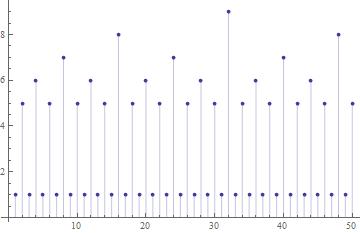}
    \caption{Graph of $v_2(15^{n}-5^{2^{8}})$ from $n$=1 to $n$=50}
    \label{fig:4.3}
\end{figure} 
\begin{figure}[H]
    \centering
    \includegraphics[width=7.5cm]{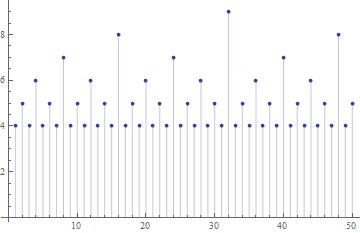}
    \caption{Graph of $v_2(n)+4$ from $n$=1 to $n$=50}
    \label{fig:4.3}
\end{figure}
Note how the behavior for even integers $n$ seems to be the same for both functions. This relation is true if and only if $m\to+\infty$ \\
Now suppose we want make
\[
v_2(a^{b}-c^{d})
\]
very large. \\
We know that $v_2(0)=+\infty$. And what if the answer is just to equate $a^{b}-c^{d}=0$ ? \cite{Rowland} I consider $b$ as a variable.
\[
a^{b}-c^{d}=0 
\]
for $b=d\frac{\log(c)}{\log(a)}$. We can expand this value to 2-adic numbers. \\ 
In fact it is possible to expand 2-adic logarithm in an infinite series of powers of base 2, more precise is the approximation, more $v_2(a^{b}-c^{d})$ grows fast. A 2-adic number is of the form:
\[
a_0+a_12+a_22^{2}+a_32^{3}+\dots
\]
\cite{p-adic}
With $0\leq a_0,a_1,a_2,a_3,\dots\leq 1$
\subsection{Example 1}
Consider the function: 
\[
v_2(15^{n}-5^{8})
\]
Bacause of what I said before:
\[
\lim_{n\to 8\frac{\log_2(5)}{\log_2(15)}\in\mathbb{Z_2}} v_2(15^{n}-5^{8}) = +\infty
\]
Using a software as pari/gp it is possible to compute the first terms of the 2-adic expansion of $8\frac{log_2(5)}{log_2(15)}$
The first eight terms of this expansion are:
\[
2 + 2^{4} + 2^{5} + 2^{8} + 2^{10} + 2^{12} + 2^{16} + 2^{17} = 202034
\]
It could be verified that:
\[
v_2(15^{202034}-5^{8}) = 26
\]
\section{Conclusions}
This result may be easily used on a calculator, to facilitate the analysis of the divisibility by 2 of this function for large values. Furthermore 2-adic valuation and in general 2-adic analysis is related to binary system, and every work related to it could have a lot of applications in informatics. p-adic numbers system perform a fundamental role in number theory, so p-adic numbers are closely related to integers. For example the extension of the 2-adic logarithm, which for more terms of the approximation gives larger values of $v_2(a^{b}-c^{d})$
\bibliographystyle{plain}
\bibliography{Bibliografia.bib}
\end{document}